\documentclass[english,letter,12pt,fleqn]{article}
\usepackage{amsmath,amssymb,amsfonts,hhline}
\usepackage{babel}
\usepackage{graphicx}
\usepackage{natbib}
\usepackage[hmargin=2.5cm,vmargin=2cm]{geometry}

\bibpunct{(}{)}{;}{a}{}{,}

\newcommand{\bs}{\boldsymbol}
\newcommand{\bo}{\mathbf}
\newcommand{\cov}{\hbox{Cov}}
\newcommand{\vect}{\hbox{vec}}
\newcommand{\ave}{\hbox{ave}}

\newtheorem{result}{Result}[section]
\newtheorem{assumption}{Assumption}[section]

\title{Nonparametric Analysis of Clustered Multivariate Data}

\author{Jaakko Nevalainen, Denis Larocque, Hannu Oja and Ilkka P{\"o}rsti
        \footnote{Jaakko Nevalainen is Professor of Biostatistics, Statistics/Department of Social Research, 20014 University of Turku, Finland (email: jaakko.nevalainen@utu.fi). Denis Larocque is Professor of Statistics, Department of Management Sciences, HEC Montr\'eal, Quebec, Canada H3T-2A7 (email: denis.larocque@hec.ca). Hannu Oja is Academy Professor, Tampere School of Public Health, 33014 University of Tampere, Finland (email: hannu.oja@uta.fi). Ilkka P{\"o}rsti is Professor of Internal Medicine, Medical School, 33014 University of Tampere, Finland. Constructive comments by the Associate Editor and two anonymous referees greatly improved the paper.
This research was supported by the Academy of Finland, The Finnish Foundation of Cardiovascular Research and the Competitive Research Funding of the Pirkanmaa Hospital District. The research work of Denis Larocque was supported by the Natural Sciences and Engineering Research Council of Canada.
}}
\date{Version: \today}
\begin{document}
\maketitle

\begin{abstract}
There has been a wide interest to extend univariate and multivariate
nonparametric procedures to clustered and hierarchical data.
Traditionally, parametric mixed models have
been used to account for the correlation structures among the dependent
observational units. In this work we extend multivariate
nonparametric procedures for one-sample and several samples location problems
to clustered data settings. The results are given for
a general score function, but with an emphasis on spatial sign and rank methods.
Mixed models notation involving design matrices for fixed and random
effects  is used throughout. The asymptotic variance formulas and limiting distributions
of the test statistics under
the null hypothesis and under a sequence of alternatives are derived,
as well as the limiting distributions for the corresponding estimates.
The approach based on a general score function also shows, for example, how
$M$-estimates behave with clustered data.
Efficiency studies demonstrate practical advantages and disadvantages
of the use of spatial sign and rank scores, and their weighted versions.
Small sample procedures based on sign change and permutation principles are
discussed. Further development of nonparametric methods for cluster correlated
data would benefit from the notation already familiar to statisticians working under
normality assumptions. Supplemental materials for the article are available online.

\vspace{10pt}\noindent\textbf{Key words and phrases:}
clustered data; mixed models; multivariate analysis; nonparametric methods.
\end{abstract}

\section{\label{section:intro}Introduction}

There has been a wide interest to extend univariate and
multivariate nonparametric procedures to clustered and hierarchical
data, which frequently arise in longitudinal studies for example.
It is well known that
unless the clustered structure is taken into account during the analysis,
the tests and confidence intervals will not maintain their prescribed levels,
leading to invalid inference.
Traditionally, parametric mixed models have been used to account for
the correlation structures among dependent observational units.
The extensions of nonparametric
methods to clustered data can roughly be divided into
univariate
\citep{RosnerGrove1999,RosnerGlynnLee2003,Rosner2006,Larocque2010,WilliamsonDattaSatten2003,Larocque2005,DattaSatten2005,Werner2007,Datta2008,LarocqueNevalainenOja2008,Kloke2009,Konietschke2009}
and multivariate approaches
\citep{Larocque2003,LarocqueNevalainenOja2007,NevalainenLarocqueOja2007,NevalainenLarocqueOja2007a,Haataja2009}.
In this paper we demonstrate that the nonparametric procedures,
which may first seem a sparse collection of tests and estimates,
can in fact be regarded as a class of score-based methods for clustered data problems. Strict
assumptions on the distribution of the random effects or the
random errors are unnecessary in this class. Our notation coincides with the one used in the mixed
models framework.

Let $\bo{Y}=(\bo{y}_1,\ldots,\bo{y}_n)'=(\bo{y}_{(1)},\ldots,\bo{y}_{(p)})$
be a sample of $p$-variate ($p>1$) random vectors with sample size $n$. The data are assumed to
be clustered throughout the paper. The cluster memberships are given by
the $n \times d$ matrix $\bo{Z}=(\bo{z}_1,\ldots,\bo{z}_n)'$:
$$\left ( \bo{Z} \right )_{ij} = \left\{
    \begin{array}{ll}
      1, & \hbox{if the $i$th observation is from cluster $j$;} \\
      0, & \hbox{otherwise.}
    \end{array}
  \right.
$$
It is useful to note that
$$\left ( \bo{Z}\bo{Z}' \right )_{ij} = \left\{
                                          \begin{array}{ll}
                                            1, & \hbox{if the $i$th and the $j$th observation are from the same cluster;} \\
                                            0, & \hbox{otherwise,}
                                          \end{array}
                                        \right.
$$
and that $\bo{Z}'\bo{Z}$ is a $d \times d$ diagonal matrix with the
cluster sizes on the diagonal, say, $m_1,\ldots,m_d$.
We also write $\bo{1}_n$ for a column $n$-vector of ones, $\vect(\bo{Y})$ for the vector
obtained by stacking the columns of $\bo{Y}$, and $\otimes$ for the Kronecker product.

A parametric linear mixed effects model for multivariate clustered data can be written as
\begin{eqnarray}
\label{parametricmodel}
\bo{Y}=\bo{Z}\bs{\alpha}+\bo{X}\bs{\beta}+\bo{E},
\end{eqnarray}
where $\bo{X}$ and $\bo{Z}$ are design matrices corresponding to the fixed effects and cluster
memberships, respectively, $\bs{\alpha}$ is a $d\times p$ random matrix of regression coefficients (random effects),
$\bs{\beta}$ is a $q\times p$ matrix of regression coefficients (fixed effects), and $\bo{E}$
is an $n \times p$ matrix of random errors. In a normality based model it is assumed that
the rows of $\bs{\alpha}$ are i.i.d. from $N_p(\bo{0},\bs{\Omega})$ and that the rows of $\bo{E}$
are i.i.d. from $N_p(\bo{0},\bs{\Sigma})$. To better illustrate the dependency structure, rewrite the
model as
\begin{eqnarray*}
\bo{Y}=\bo{X}\bs{\beta}+\bo{E}_* \hbox{, where } \vect(\bo{E}_*')\sim N_{np}(\bo{0},\bo{I}_n\otimes\bs{\Sigma}+\bo{Z}\bo{Z}'\otimes \bs{\Omega}).
\end{eqnarray*}
If $\bo{E}_*=(\bs{\epsilon}_1,\ldots,\bs{\epsilon}_n)'$ this means that
\begin{itemize}
\item[(P1)] $\bs{\epsilon}_i \sim N_p(\bo{0},\bs{\Sigma}+\bs{\Omega})$ for all $i=1,\ldots,n$.
\item[(P2)] If $(\bo{Z}\bo{Z}')_{ij}=1$ then
$$\vect(\bs{\epsilon}_i,\bs{\epsilon}_j)\sim N_{2p}\left(\bo{0}, \left(
                                                                   \begin{array}{cc}
                                                                     \bs{\Sigma}+\bs{\Omega} & \bs{\Omega} \\
                                                                     \bs{\Omega} & \bs{\Sigma}+\bs{\Omega} \\
                                                                   \end{array}
                                                                 \right)
 \right ).$$
\item[(P3)] If $(\bo{Z}\bo{Z}')_{ij}=0$ then $\bs{\epsilon}_i$ and $\bs{\epsilon}_j$ are independent.
\end{itemize}
The problem with the parametric linear model are the strict distributional assumptions, which
can be completely unrealistic. Its use outside of the assumed model can lead into inefficient or even invalid statistical inference.  In section \ref{section:score_approach} of this paper we introduce
an alternative semiparametric linear model to analyze clustered data under relaxed assumptions.
Section \ref{section:onesample}
treats the multivariate one-sample location problem, both from the point of
view of testing and estimation. The section has some review character in it but we now present
the results for a general score function rather than focusing on a specific score.
In section \ref{section:multisample} the treatment
is expanded to the multivariate several samples location problem, which has not appeared
previously. Efficiency studies
demonstrate the practical advantages and disadvantages of the proposed
two-sample weighted nonparametric tests compared to the classical methods based
on the sample mean and the sample covariance matrix in section \ref{section:simlations}.
Finally, two data sets are analyzed with spatial sign and rank methods (provided as a supplemental file).

\section{\label{section:score_approach}A Semiparametric Linear Model}

Suppose that the data is clustered and we wish to analyze it by the linear model
\begin{eqnarray*}
\bo{Y}=\bo{X}\bs{\beta}+\bo{E}
\hbox{, where } \bo{E}=(\bs{\epsilon}_1,\ldots,\bs{\epsilon}_n)'.
\end{eqnarray*}
However, we wish to avoid assumptions of
normality, or any other parametric distribution, on the random errors. Instead we assume the following.
\begin{assumption}Distributional assumptions.
\begin{itemize}
  \item[(D1)] The marginal distributions are identical:
  $\bs{\epsilon}_i \sim \bs{\epsilon}_j$ for all $i, j=1,\ldots,n$.
  \item[(D2)] The pairwise joint distributions are identical: $\vect(\bs{\epsilon}_i,\bs{\epsilon}_j) \sim \vect(\bs{\epsilon}_{i'}, \bs{\epsilon}_{j'})$ for all $i \ne j$ and $i' \ne j'$ and $(\bo{Z}\bo{Z}')_{ij}=(\bo{Z}\bo{Z}')_{i'j'}$.
  \item[(D3)] The clusters are independent: if $(\bo{Z}\bo{Z}')_{ij}=\bo{0}$ then $\bs{\epsilon}_i$ and $\bs{\epsilon}_j$ are independent.
\end{itemize}
\end{assumption}
Assumptions (D1) and (D2) fix the location, variance and covariance structure of the model. Assumption (D3) is a natural and standard
presumption, and also needed
for finding the limiting distribution. Compared to the assumptions of the parametric linear model (P1)-(P3), these conditions are
not restrictive.

To work under these assumptions it is often advantageous to use transformed observations instead of the original ones.
Our approach is first to apply a vector-valued score function $\bo{T}(\bo{y})$
to the data points (or to the centered observations as discussed later). To prove the asymptotic results,
it is sufficient that the score function has the following properties.
\begin{assumption}
Sufficient conditions on the score function.
\begin{itemize}
\item[(S1)] The score function satisfies
$$\|\bo{T}(\bo{y}+c\bs{\delta}) - \bo{T}(\bo{y}) -  \dot{\bo{T}}_{\bs{\delta}}(\bo{y})c\| \leq R_{\bs{\delta}}(\bo{y})c^{1+\zeta}$$
for some $\zeta>0$, and $\exists K >0$ and $M>0$ such that $E \left [ \sup_{\|\bs{\delta}\|\leq K} R_{\bs{\delta}}(\bs{\epsilon}_i) \right ] < M.$
\item[(S2)] $E \left ( \|\bo{T}(\bs{\epsilon})\|^\nu \right )<\infty$  for some $\nu>2$.
\end{itemize}
In the one-sample case, one additional condition is needed.
\begin{itemize}
\item[(S3)] The score function is odd: $\bo{T}(-\bo{y})=-\bo{T}(\bo{y})$.
\end{itemize}
\end{assumption}
Note that if $\nabla \bo{T}(\bo y)$ is the matrix of partial derivatives
$\partial \bo T_i(\bo y)/\partial \bo y_j$, then
$\dot{\bo{T}}_{\bs{\delta}}(\bo{y})=\nabla \bo{T}(\bo y) \bs{\delta}$ in (S1).

\begin{figure}[h!]
\includegraphics[width=\textwidth]{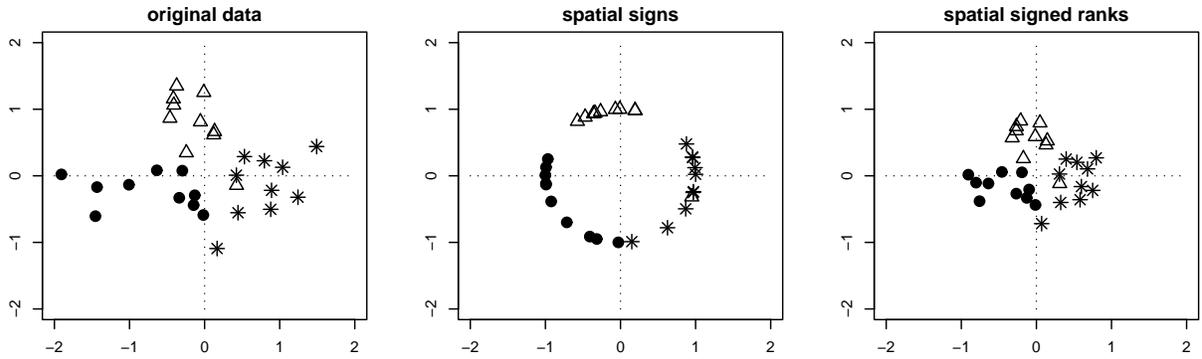}
\caption{\label{figure:clustereddata}Clustered data from a bivariate
spherical normal distribution with intracluster correlation of 2/3.
The three clusters are indicated by different symbols.  The spatial
signs of the observations from the same cluster tend to lie on the
same regions of the unit circle, whereas the signed ranks also preserve
the shape of the data cloud, and the cluster structure remains
clearly visible.}
\end{figure}

Tests and estimates can then be constructed on the transformed data.
By choosing the score function well, tests and estimates can achieve desired properties for the problem
at hand, like robustness against outliers, or improved efficiency for
heavy-tailed distributions.  Some examples of clever choices of scores are signs and ranks
commonly used in nonparametric statistics, optimal scores from maximum likelihood theory,
or Huber's score often applied in robust statistics.
Thus, the proposed approach works for a general score function, and suggests how the
tests and estimates could be constructed. Some authors have expressed interest towards this type of approach in univariate testing \citep{Jin2003,Huang2008},
but these tests are currently designed for independent observations only.
The main motivation for the present paper is, however, multivariate sign and rank methods
resulting from taking either the \textit{spatial sign}
$$\bo{S}(\bo{y}_i)=\|\bo{y}_i\|^{-1}\bo{y}_i,$$ the centered \textit{spatial
rank}
$$\bo{R}_n(\bo{y}_i)=\frac{1}{n}\sum_{j=1}^n
\|\bo{y}_i-\bo{y}_j\|^{-1}(\bo{y}_i-\bo{y}_j) =
\frac{1}{n}\sum_{j=1}^n \bo{S}(\bo{y}_i-\bo{y}_j),$$ or the centered
\textit{spatial signed rank}
\begin{eqnarray*}
\bo{Q}_n(\bo{y}_i) & = & \frac{1}{2n}\sum_{j=1}^n \left [
\|\bo{y}_i-\bo{y}_j\|^{-1}(\bo{y}_i-\bo{y}_j) +
\|\bo{y}_i+\bo{y}_j\|^{-1}(\bo{y}_i+\bo{y}_j) \right  ] \\ & = &
\frac{1}{2} \left [ \bo{R}_n(\bo{y}_i) - \bo{R}_n(-\bo{y}_i)\right
].
\end{eqnarray*} By convention, $\|\bo{0}\|^{-1}\bo{0}=\bo{0}$. Figure
\ref{figure:clustereddata} illustrates how the transformations
preserve the clustering structure. The methods based on spatial
signs and ranks are more robust, more efficient for heavy-tailed
distributions than normal theory based methods
\citep{MottonenOjaTienari1997}, and do not require assumptions on
the existence of moments of $\bo{y}_i$. Spatial sign and rank methods have
been criticized for their lack of affine invariance and equivariance properties,
but this problem can be overcome by a modified transformation-retransformation procedure (section \ref{section:conclusions}).
The results of the paper have been written having these three score functions
in mind but they hold more generally.

Conditions (S1)--(S3) hold for the identity score $\bo{T}(\bo{y})=\bo{y}$ if
$E \left ( \|\bo{y}\|^\nu \right )<\infty$  for some $\nu>2$.
For the spatial sign score $\bo{S}(\bo{y})$, (S2) and (S3) are trivially true. To verify (S1), one can first show that
$$\big\| \bo{S}(\bo{y}+c\bs{\delta})-\bo{S}(\bo{y})\frac{1}{\|\bo{y}\|}\left [\bo{I}_p-\bo{S}(\bo{y})\bo{S}(\bo{y})' \right ]\bs{\delta}c {\big\|}  \leq B \frac{\|\bs{\delta}\|^{1+\zeta}}{\|\bo{y}\|^{1+\zeta}}c^{1+\zeta}, 0 < \zeta <1.$$
\citep{Arcones1998,BaiChen1990}. (S1) then follows if $\bs{\epsilon}_i$ has a bounded density.
The conditions can be established similarly for the spatial rank score.

\section{\label{section:onesample}One-Sample Case}

Assume that
$$\bo{Y} = \bo{1}_n\bs{\mu}' + \bo{E},$$
where rows of $\bo{E}$ satisfy assumptions (D1)-(D3).
In the one-sample case it is natural to transform the data set $\bo{Y}=(\bo{y}_1,\ldots,\bo{y}_n)' \rightarrow \bo{T}=(\bo{T}(\bo{y}_1),\ldots,\bo{T}(\bo{y}_n))'$ using an odd score function $\bo{T}(\bo{y})$.
We wish to test the null hypothesis $H_0: \bs{\mu}=\bo{0}$ without loss of generality, where the
location parameter $\bs{\mu}$ satisfies $E\left ( \bo{T}(\bo{y}_i-\bs{\mu})\right )=\bo{0}$. Thus, its
interpretation depends on the choice of the score. However, if the distribution is symmetric, $\bo{y}_i -\bs{\mu} \sim \bs{\mu}-\bo{y}_i$,
all tests test the same null hypothesis and the corresponding estimates estimate the same population parameter with different statistical
properties.

Let $f$ be the density of $\bs{\epsilon}_i$ and $\bo{L}(\bo{y}_i)$ be the optimal score function, the gradient vector of $\log f(\bo{y}-\bs{\mu})$
with respect to $\bs{\mu}$ at the origin. Note also that $E \left (\bo{T}(\bo{y}_i) \right )=\bo{0}$ if $H_0$ is true.
Define
$$\bo{A} = E \left ( \bo{T}(\bs{\epsilon}_i) \bo{L}(\bs{\epsilon}_i) ' \right ) \hbox{ and }
\bo{B} = E \left ( \bo{T}(\bs{\epsilon}_i) \bo{T}(\bs{\epsilon}_i)' \right )$$
and
$$\bo{C} = E \left ( \bo{T}(\bs{\epsilon}_i) \bo{T}(\bs{\epsilon}_j)' \right ) \hbox{ where } i \ne j \hbox{ satisfy } (\bo{Z}\bo{Z}')_{ij}=1.$$
The covariance structure of the weighted scores $\bo{W}\bo{T}$ is then
$$\cov \left ( \vect(\bo{T}'\bo{W}) \right ) = \bo{W}^2 \otimes \bo{B} + \left (\bo{W} ( \bo{Z}\bo{Z}' - \bo{I}_n ) \bo{W}\right )
\otimes \bo{C},$$
where $\bo{W}=diag(\bo{w})$ is a $n \times n$ diagonal matrix with a non-negative weight associated with the $i$th observation as the $i$th diagonal element. For the sampling design and the weights, it is assumed that $\bo{1}_n' \bo{W} \bo{1}_n = n$ and that there exist constants $D_{B}$ and $D_C$ such that
$$\frac{1}{n} \bo{1}_n' \bo{W}^2 \bo{1}_n \rightarrow D_{B} \hbox{ and }
\frac{1}{n} \bo{1}_n' \bo{W} \left ( \bo{Z}\bo{Z}' - \bo{I}_n \right
) \bo{W} \bo{1}_n \rightarrow D_{C}$$ as $d$ tends to infinity.

But how should the weights be chosen? It is natural that the members of the same cluster receive the same weight.
Furthermore, it can be shown that if $\bo{C}=\rho \bo{B}$ and the covariance matrix has the structure
$$\cov(\vect(\bo{T}'))=\bo{I}_n \otimes \bo{B}+(\bo{Z}\bo{Z}'-\bo{I}_n) \otimes \bo{C}
= \bs{\Sigma} \otimes \bo{B}, \hbox{ where }\bs{\Sigma} = \bo{I}_n+\rho(\bo{Z}\bo{Z}'-\bo{I}_n),$$
the optimal weights are given by  $\bo{w} = \kappa \bs{\Sigma}^{-1}\bo{1}_n$. The weights in the $i$th cluster are then proportional
to $\left [ 1+(m_i-1) \rho \right ]^{-1}$ \citep{LarocqueNevalainenOja2007}.
These weights are optimal in the sense that they maximize the Pitman efficiency of a test based on $\bo{W}\bo{T}$. Here $\kappa$ is the Lagrange multiplier chosen so that the constraint $\bo{w}' \bo{1}_n=n$ is satisfied.

\subsection{Testing}

The test statistic for testing $H_0: \bs{\mu} = \bo{0}$ versus $H_1:
\bs{\mu} \neq \bo{0}$ is the weighted average of the scores
$$\frac{1}{n}\bo{T}'\bo{W}\bo{1}_n  = \frac{1}{n}\left ( \bo{1}_n' \otimes \bo{I}_p \right ) \vect \left (\bo{T}'\bo{W} \right ).$$
Then we have:

\begin{result} Under the null hypothesis $H_0: \bs{\mu} = \bo{0}$,  as $d$ tends to
infinity
      $$Q^2 = \frac{1}{n}\bo{1}_n' \bo{W} \bo{T} \left ( \frac{1}{n} \bo{T}' \bo{W} \bo{Z} \bo{Z}' \bo{W} \bo{T} \right )^{-1} \bo{T}' \bo{W} \bo{1}_n \stackrel{d}\rightarrow \chi^2_p,$$
      where $n^{-1} \bo{T}' \bo{W} \bo{Z} \bo{Z}' \bo{W} \bo{T}$ is a consistent estimate
      of $D_{B}\, \bo{B}+D_{C}\, \bo{C}$, the asymptotic covariance matrix of
      $\frac{1}{\sqrt{n}}  \bo{T}'\bo{W}\bo{1}_n$.
\end{result}

The result follows from the generalization of the central limit theorem given as a corollary in \citet[p. 30]{Serfling1980},
by noting that the cluster sums are independent but not identically distributed random variables.

For symmetric distributions, small sample (meaning here a small $d$) $p$-values can be based on
the sign change principle. Under the null,
$$\bo{T}'\bo{W}\bo{1}_n = \bo{T}'\bo{W} \bo{Z} \bo{1}_d \sim \bo{T}'\bo{W} \bo{Z} \bo{J}_d \bo{1}_d,$$ where
$\bo{J}_d$ is a $d \times d$ diagonal sign-change matrix, with $\pm 1$ as diagonal entries, and changing all the signs within a cluster at the same time. The covariance matrix estimate is invariant under these sign changes, and the test statistic becomes
$$Q^2_J = \frac{1}{n}\bo{1}_d' \bo{J}_d \bo{Z}' \bo{W} \bo{T} \left ( \frac{1}{n} \bo{T}' \bo{W} \bo{Z} \bo{Z}' \bo{W} \bo{T} \right )^{-1} \bo{T}' \bo{W} \bo{Z} \bo{J}_d \bo{1}_d.$$
Estimation of $E_J \left [ I( Q^2_J \ge Q^2) | \bo{Y} \right ]$ over \textit{e.g.} 1000 equiprobable random allocations of signs gives an
estimated $p$-value for the conditionally distribution-free test.

\begin{result} The limiting distribution under the
sequence of alternatives $H_n: \bs{\mu} = n^{-1/2}\bs{\delta}$ is a
non-central chi-square
      $$Q^2 \stackrel{d}\rightarrow \chi^2_p \left (\bs{\delta}' \bo{A} \left( D_{B}\, \bo{B}+D_{C}\, \bo{C} \right ) ^{-1} \bo{A} \bs{\delta} \right )$$
as $d\rightarrow \infty$.
\end{result}
This result allows the computation of asymptotic relative efficiencies \citep{LarocqueNevalainenOja2007,Haataja2009}.

\subsection{\label{section:onesampleestimation}Estimation}

The companion estimate of location is determined by the estimating
equation
$$\bo{1}_n' \bo{W} \hat{\bo{T}} = \bo{0},$$
where $\hat{\bo{T}} = (\bo{T}(\bo{y}_1- \hat{\bs{\mu}}),\ldots,
\bo{T}(\bo{y}_n- \hat{\bs{\mu}}))'$. Thus, the solution is a
location estimate with the property that the weighted scores of
shifted observations add up to zero. If (S1) holds and $\sqrt{n} \left
(\hat{\bs{\mu}}-\bs{\mu} \right )=\bo{O}_p(1)$, the Bahadur-type representation
$$\sqrt{n} \left
(\hat{\bs{\mu}}-\bs{\mu} \right ) = \frac{1}{\sqrt{n}} \bo{A}^{-1}
\bo{T}' \bo{W} \bo{1}_n + \bo{o}_p(1)$$
shows the relationship between the test and the estimate. The connection has been established in detail
for the spatial sign score \citep{NevalainenLarocqueOja2007,NevalainenLarocqueOja2007a}
with clustered data, and for the spatial signed-rank score \citep{Chaudhuri1992}
with independent observations.
Now, the asymptotic distributions of the weighted spatial median or the
weighted spatial Hodges-Lehmann estimate for example, are trivial:
\begin{result}
The limiting distribution of the score-based estimate is
$$\sqrt{n} \left (\hat{\bs{\mu}}-\bs{\mu} \right )
\stackrel{d}\rightarrow N_p \left (\bo{0}, \bo{A}^{-1} \left ( D_B \bo{B} + D_C \bo{C} \right ) \bo{A}^{-1} \right )$$
as $d \rightarrow \infty$.
\end{result}
The estimation of the covariance matrix can be based on the residuals $\hat{\bs{\epsilon}}_i=\bo{y}_i-\hat{\bs{\mu}}$.
The matrix $D_B \bo{B} + D_C \bo{C}$ can be consistently estimated by $n^{-1} \hat{\bo{T}}' \bo{W} \bo{Z} \bo{Z}' \bo{W} \hat{\bo{T}}$.
For the weighted spatial median and the weighted Hodges-Lehmann estimate the  $\bo{A}$ matrix is estimated by
$$\ave_i \left [ \bo{A}(\hat{\bs{\epsilon}}_i) \right ] \mbox{ and }
\ave_{ij} \left [ \bo{A} \left ( \frac{1}{2} (\hat{\bs{\epsilon}}_i+\hat{\bs{\epsilon}}_j)\right ) \right ]\hbox{, where }
\bo{A}(\bs{\epsilon}_i)=\frac{1}{\|\bs{\epsilon}_i\|} \left ( \bo{I}_p-\frac{\bs{\epsilon}_i\bs{\epsilon}_i'}{\|\bs{\epsilon}_i\|^2}\right ),$$
and the second average is over the pairs with $(\bo{Z}\bo{Z}')_{ij}=0$.
Estimation of precision of the estimates in small samples could potentially be based on
bootstrap procedures for clustered data \citep{Field2007}, but more practical experience is needed on this approach.

\section{\label{section:multisample}Several Samples Case}

Suppose now that the data consist of $(\bo{X},\bo{Z},\bo{Y})$, where $\bo{Y}$ and $\bo{Z}$ are matrices of response vectors and the cluster memberships, respectively, in the same way as in section \ref{section:intro}, and that $\bo{X}=(\bo{x}_1,\ldots,\bo{x}_n)'$ is an $n \times c$ matrix indicating group or sample membership such that
$$\left ( \bo{X} \right )_{ij} = \left\{
    \begin{array}{ll}
      1, & \hbox{if the $i$th observation is from group $j$;} \\
      0, & \hbox{otherwise.}
    \end{array}
  \right.
$$
Again,
$$\left ( \bo{X}\bo{X}' \right )_{ij} = \left\{
                                          \begin{array}{ll}
                                            1, & \hbox{if the $i$th and the $j$th observation are from the same group;} \\
                                            0, & \hbox{otherwise,}
                                          \end{array}
                                        \right.
$$
and that $\bo{X}'\bo{X}$ is a $c \times c$ diagonal matrix with the group sizes
on the diagonal, say, $n_1,\ldots,n_c$, again assumed fixed by the design. Recall the cluster membership
matrix $\bo{Z}$ and the weight matrix $\bo{W}$ from the previous section. Now $\bo{X}'\bo{Z}$ is a frequency
table fixing the design.

Write the model as
$$\bo{Y} = \bo{1}_n\bs{\mu}' + \bo{X}\bs{\Delta}' + \bo{E},$$
where $\bs{\mu}$ is the overall location center (\textit{e.g.} grand mean), $\bs{\Delta}$
is a $p \times c$ contrast matrix representing the treatment effects or the
deviations from that location center, and the rows of $\bo{E}$ satisfy (D1)-(D3).
The parameters in the model depend on the choice of the score, population, and the design.
For now we are interested in the treatment effects only: the goal is to confront the hypotheses
$H_0: \bs{\Delta} = \bo{0} \hbox{ vs. } H_1: \bs{\Delta} \ne \bo{0}.$

 In the test
construction we need estimated (or centered) scores $\hat{\bo{T}}$
chosen to satisfy $\bo{1}_n' \bo{W} \hat{\bo{T}} = \bo{0}$.
Different scores require different inner centering to fulfil this property. For example, for
the identity or the spatial sign score, the estimated scores and the theoretical scores are
$\hat{\bo{T}}_i=\bo{T}\left ( \bo{y}_i - \hat{\bs{\mu}} \right ) \hbox{ and } {\bo{T}}_i=\bo{T}\left ( \bo{y}_i - {\bs{\mu}} \right )$,
where $\bs{\mu}$ is estimated from the whole sample assuming absence of treatment effects.
As will be seen in the next section, the test statistic can be expressed in two asymptotically equivalent forms, but only if the weighted estimate of location is based on the same score as the test.
For these two scores, the estimate of $\bs{\mu}$ should be the weighted
sample mean (identity score) or the weighted spatial median (spatial sign score). Rank scores are automatically centered but the
weighted ranks are not. Thus, in the case of ranks we write
$$\hat{\bo{T}}_i = n^{-1} \sum_{j=1}^n w_j \bo{S}(\bo{y}_i - \bo{y}_j) \hbox{ and }  {\bo{T}}_i = E\left (\bo{S}(\bo{y}_i-\bo{y}_j) \right ),$$
where the expectation is taken over $\bo{y}_j$ with $(\bo{Z}\bo{Z}')_{ij}=0$.
Note that $\sum_{i=1}^n w_i \hat{\bo{T}}_i = \bo{0}$ and $E({\bo{T}}_i)=\bo{0}$. Collect these and let
$\hat{\bo{T}} = \left ( \hat{\bo{T}}_1,\ldots,\hat{\bo{T}}_n \right )'$   and  $\bo{T}=\left ( \bo{T}_1,\ldots,\bo{T}_n \right )'$
 denote the matrices of  estimated and theoretical scores, respectively. Under the null hypothesis, $E \left ( \vect (\bo{T}'\bo{W}) \right ) = \bo{0}$ and the covariance structure is given by
$$\cov \left ( \vect  \left ( \bo{T}'\bo{W} \right ) \right ) = \bo{W}^2 \otimes \bo{B} + \left (\bo{W} \left ( \bo{Z} \bo{Z}'-\bo{I}_n \right ) \bo{W} \right )\otimes \bo{C}.$$

Optimal weights are obtained for the two-sample problem in a similar way as for the one-sample problem: if again $\bo{C}=\rho\bo{B}$
and $\bo{G}=diag(\bo{x}_{(1)}-\bo{x}_{(2)})$ the optimal weights are
$\bo{w}=\bs{\Sigma}^{-1}\left(\kappa_1\bo{x}_{(1)} + \kappa_2\bo{1}_n \right ),$
where $\bs{\Sigma}=(1-\rho)\bo{I}_n+\rho\bo{G}\bo{Z}\bo{Z}'\bo{G}$
with Lagrange multipliers $\kappa_1$ and $\kappa_2$ chosen such that  $\bo{w}'\bo{x}_{(1)} = n_1$ and $\bo{w}' \bo{1}_n=n$.

\subsection{\label{subsection:severalsamplestesting}Testing}

The test to confront the hypotheses $H_0: \bs{\Delta}=\bo{0}$ and $H_1: \bs{\Delta} \neq \bo{0}$ is based on the $p(c-1)$-vector
$$\widetilde{\bo{T}} = \frac{1}{\sqrt{n}} \left ( \bo{H}'\bo{X}' \otimes \bo{I}_p \right ) \vect \left ( \hat{\bo{T}}' \bo{W} \right ),$$
where $\bo{H}$ is a  $c \times (c-1)$ matrix obtained from the
identity matrix by dropping its $c$th column. Under the null hypothesis,
$$\widetilde{\bo{T}} = \frac{1}{\sqrt{n}} \left ( \bo{H}'\bo{X}_0' \otimes \bo{I}_p \right ) \vect (\bo{T}'\bo{W}) + \bo{o}_p(1),$$
where $\bo{X}_0 = \left
( \bo{I}_n - n^{-1} \bo{1}_n \bo{1}_n' \bo{W} \right ) \bo{X}$ is a
centered version of the design matrix, which has enabled us to
replace the estimated scores $\hat{\bo{T}}$ by the theoretical
scores $\bo{T}$.

Under the null hypothesis $E \left ( \frac{1}{\sqrt{n}} \left ( \bo{H}'\bo{X}_0' \otimes \bo{I}_p \right ) \vect (\bo{T}'\bo{W})\right ) =
\bo{0}$ and correspondingly its nonsingular covariance matrix is
$$\frac{1}{n} \left ( \bo{H}'\bo{X}_0' \bo{W}^2 \bo{X}_0 \bo{H} \right ) \otimes \bo{B} + \frac{1}{n} \left ( \bo{H}' \bo{X}_0'\bo{W} \left ( \bo{Z} \bo{Z}' - \bo{I}_n \right ) \bo{W} \bo{X}_0 \bo{H} \right ) \otimes \bo{C}.$$
A necessary assumption is that the matrices
$$\frac{1}{n} \left ( \bo{X}_0' \bo{W}^2 \bo{X}_0 \right ),
\frac{1}{n} \left ( \bo{X}_0'\bo{W} \left ( \bo{Z} \bo{Z}' -
\bo{I}_n \right ) \bo{W} \bo{X}_0 \right ) \hbox{ and } \frac{1}{n} \bo{X}'\bo{W}\bo{X}$$ converge to finite
matrix-valued limits, $\bo{D}_B$, $\bo{D}_C$ and $\bs{\Lambda}$, say. Note the
similarity of the requirement to the one-sample case: the limits now
need to exist groupwise. Diagonal elements of the matrix $\bs{\Lambda}=diag(\bs{\lambda})$ satisfy
$0 < \lambda_1,\ldots,\lambda_c<1$.

\begin{result}  Under the null hypothesis $H_0: \bs{\Delta} = \bo{0}$
      \begin{eqnarray*}
Q^2 & = & \widetilde{\bo{T}}'\left [ \left (  \frac{1}{n} \bo{H}'  \left ( \bo{X}_0'  \bo{W}^2 \bo{X}_0 \right ) \bo{H} \right ) \otimes \left ( \frac{1}{n} \hat{\bo{T}}'\hat{\bo{T}} \right ) \right . \\ & & \left . + \left (\frac{1}{n} \bo{H}'  \left ( \bo{X}_0' \bo{W} \left ( \bo{Z} \bo{Z}' -
\bo{I}_n \right ) \bo{W} \bo{X}_0 \right ) \bo{H} \right )\otimes \left (\frac{1}{k} \hat{\bo{T}}' \left ( \bo{Z}\bo{Z}' - \bo{I}_n \right )  \hat{\bo{T}} \right ) \right ]^{-1} \widetilde{\bo{T}} \stackrel{d}\rightarrow \chi^2_{p(c-1)},
\end{eqnarray*}
as $d \rightarrow \infty$, where $k=\bo{1}_n'\left( \bo{Z}\bo{Z}'-\bo{I}_n \right) \bo{1}_n$ and
the part inside the brackets $[\cdot]$ is a consistent estimate of $(\bo{H}'\bo{D}_B \bo{H}) \otimes \bo{B} +(\bo{H}'\bo{D}_C\bo{H}) \otimes \bo{C}$.
\end{result}
A conditionally distribution-free permutation test is constructed as follows. Let $\bo{P}$ be an ``acceptable"
$n \times n$
permutation matrix obtained by permuting the rows or columns of an identity matrix, uniform among acceptable permutations. The $p$-value of
the permutation test is the estimate of
$$E_P \left [ I \left ( Q^2 (\bo{P}\bo{X},\bo{Z},\bo{Y}) \ge Q^2 (\bo{X},\bo{Z},\bo{Y}) \right ) | \bo{Y} \right ].$$
The proper way to permute clustered data in an acceptable way depends on the design.
Clearly, permutations which do not change the distribution of the test statistic when the null hypothesis is true are guaranteed to provide a valid and distribution-free test.
We say that two designs $(\bo{X}_1,\bo{Z}_1)$ and $(\bo{X}_2,\bo{Z}_2)$ are \emph{equivalent in structure} if there exist permutation matrices $\bo{P}_c$ and $\bo{P}_d$ of dimensions $c\times c$ and $d \times d$, respectively, such that $$\bo{P}_c \bo{X}_1' \bo{Z}_1 \bo{P}_d = \bo{X}_2' \bo{Z}_2.$$
To ensure that the null distribution of the test
statistic is invariant under permutations, the general condition on the permutation matrix is that  $(\bo{P}\bo{X},\bo{Z})$ and $(\bo{X},\bo{Z})$ are equivalent in structure. In controlled trials the approach should also follow the randomization scheme.
One can distinguish between three common designs:

\bigskip\noindent \textbf{Design A.} A permutation fulfilling the general condition is natural for observational studies, where the researcher has no control over the group memberships.

\bigskip\noindent \textbf{Design B.} Randomization of individuals inside the clusters. The permutation of the treatment assignments should be performed only within the clusters. More formally, so that the permutation matrices $\bo{P}$ satisfy $\bo{P}\bo{Z}=\bo{Z}$. The condition implies the general condition.

\bigskip\noindent \textbf{Design C.} Randomization of clusters. The permutation should then maintain the members of the same cluster within the same treatment, and the permutation matrices $\bo{P}$ should satisfy the general condition and
    $\bo{P}\bo{Z}\bo{Z}'\bo{P}'=\bo{Z}\bo{Z}'.$ This allows the exchange of treatments between clusters of the same size only, and the permutation may not be very rich.

Consider next the sequence of alternatives $H_n:
\bs{\Delta} = n^{-1/2} \bs{\Delta}_0$, where $\bs{\Delta}_0=(\bs{\delta}_1,\ldots,\bs{\delta}_c)$ is
chosen in such a way that $\bs{\Delta}_0 \bs{\lambda} =  \bo{0}$.  This fixes the location parameter $\bs{\mu}$ for asymptotic studies.
Under $H_n$,
$$E \left ( \widetilde{\bo{T}} \right ) = \frac{1}{n}\left ( (\bo{H}' \bo{X}_0' \bo{W}) \otimes \bo{I}_p \right ) \vect \left ( \bo{A} \bs{\Delta}_0 \bo{X}'\right ) = \frac{1}{n} \left ( (\bo{H}'\bo{X}'\bo{W} \bo{X} ) \otimes \bo{A} \right ) \vect \left ( \bs{\Delta}_0\right ) +\bo{o}(1).$$
Therefore we have the following result.
\begin{result}
Under $H_n$, the limiting distribution of $Q^2$ is a noncentral chi-square with $p(c-1)$ degrees of freedom and noncentrality parameter
\begin{eqnarray*}
\vect(\bs{\Delta}_0 \bs{\Lambda} \bo{H})' \left [ (\bo{H}'\bo{D}_B \bo{H}) \otimes \left ( \bo{A}^{-1}\bo{B}\bo{A}^{-1} \right ) +  (\bo{H}'\bo{D}_C\bo{H}) \otimes \left ( \bo{A}^{-1}\bo{C}\bo{A}^{-1} \right ) \right ]^{-1}\vect(\bs{\Delta}_0 \bs{\Lambda} \bo{H}),
\end{eqnarray*}
as $d \rightarrow \infty$.
\end{result}
Alternatively, the noncentrality parameter can be written as
\begin{eqnarray*}
\vect(\bs{\Delta}_0)' \left [ \bo{D}_B \otimes \left ( \bo{A}^{-1}\bo{B}\bo{A}^{-1} \right ) + \bo{D}_C \otimes \left ( \bo{A}^{-1}\bo{C}\bo{A}^{-1} \right )\right ]^+ \vect(\bs{\Delta}_0)\\
=\vect(\bo{A}\bs{\Delta}_0)' \left [ \bo{D}_B \otimes \bo{B} + \bo{D}_C \otimes \bo{C}\right ]^+ \vect(\bo{A}\bs{\Delta}_0),
\end{eqnarray*}
which is of the same form as in the one-sample case.

\subsection{Estimation}

Until now we have parametrizised the model with $\bs{\mu}$ and $\bs{\Delta}$,
which depend not only on the underlying population but also on the design. Let us now reparametrize the model by $\bs{\beta} = (\bs{\mu}_1,\ldots,\bs{\mu}_c)' = \bo{1}_c\bs{\mu}' + \bs{\Delta}'$, and obtain the model
$$\bo{Y} = \bo{X}\bs{\beta} + \bo{E}.$$
The weighted estimates of the group centers
$\bs{\mu}_1,\ldots,\bs{\mu}_c$ are found via solving the $c$
estimating equations $\bo{X}' \bo{W} \hat{\bo{T}} = \bo{0},$ where
now $\hat{\bo{T}}_i=\bo{T} ( \bo{y}_{i} - \hat{\bs{\beta}}'\bo{x}_i  )$, and $\hat{\bo{T}} = \left (
\hat{\bo{T}}_1,\ldots,\hat{\bo{T}}_n \right )'$. In essence, this is simply
a one-sample estimation problem repeated $c$ times (section  \ref{section:onesampleestimation}).

Estimation of group differences is a little more subtle issue. Due to clustering, the
observations are correlated, and so are the estimates. Write $\bs{\theta}_{ij} =
{\bs{\mu}}_j - {\bs{\mu}}_i$, $i,j=1,\ldots,c$. With the identity, spatial sign and rank scores, the problem is reduced to computation
of (i) the difference of the mean vectors, (ii) the difference of the spatial medians, or (iii)
the two-sample Hodges-Lehmann estimate \citep{HodgesLehmann1963,MottonenOja1995}.
Again, under sufficient conditions, the connection between the estimate and scores is
\begin{eqnarray*}
\sqrt{n} (\hat{\bs{\theta}}_{ij}-\bs{\theta}_{ij}) & = &
\sqrt{n} \bo{A}^{-1} \left (\frac{1}{n_j} \bo{T}'\bo{W} \bo{x}_{(j)} - \frac{1}{n_i} \bo{T}'\bo{W} \bo{x}_{(i)} \right ) + \bo{o}_p(1)\\
&=& \sqrt{n} \bo{A}^{-1} \bo{T}'\bo{w}_{ij}^* + \bo{o}_p(1),
\end{eqnarray*}
where $n_i=\bo{x}_{(i)}'\bo{x}_{(i)}$ is the group size.
Standard theory yields:
\begin{result}
Under general assumptions, the limiting distribution of
$\sqrt{n} (\hat{\bs{\theta}}_{ij}-\bs{\theta}_{ij})$ is a $p$-variate normal distribution with expectation zero and covariance matrix
\begin{eqnarray}
\label{eqn:covmat_of_pairwise_diff}
\bo{A}^{-1} \left ( \gamma_B  \bo{B} + \gamma_C \bo{C} \right ) \bo{A}^{-1},
\end{eqnarray}
as $d \rightarrow \infty$ and where
\begin{eqnarray*}
\gamma_B & = & \lim_{d \rightarrow \infty} \left [ \frac{n}{n_i} \frac{\bo{x}_{(i)}' \bo{W}^2 \bo{x}_{(i)}}{n_i}  +  \frac{n}{n_j} \frac{\bo{x}_{(j)}' \bo{W}^2 \bo{x}_{(j)}}{n_j} \right ]\\
\gamma_C & = & \lim_{d \rightarrow \infty} \left [ \frac{n}{n_i} \frac{\bo{x}_{(i)}' \bo{W}(\bo{Z}\bo{Z}'-\bo{I}_n) \bo{W} \bo{x}_{(i)}}{n_i}  + \frac{n}{n_j} \frac{\bo{x}_{(j)}' \bo{W}(\bo{Z}\bo{Z}'-\bo{I}_n) \bo{W} \bo{x}_{(j)} }{n_j } \right . \\ & & \left .- 2 \frac{n}{\sqrt{n_i n_j}} \frac{\bo{x}_{(i)}' \bo{W} \bo{Z} \bo{Z}' \bo{W} \bo{x}_{(j)}}{\sqrt{n_i n_j}} \right ].
\end{eqnarray*}
\end{result}
This covariance breakdown shows how the intracluster dependency affects the
covariance structure via members of the same cluster receiving the
same and different treatments. If all members of the cluster belong to the same group,
the last part of $\gamma_C$ disappears, and the total variance can be seriously underestimated
if the clustering is ignored. The opposite may happen when treatments are assigned within the clusters.

In practice, the limiting constants $\gamma_B$ and $\gamma_C$ can be replaced by their empirical counterparts.
The estimation of the matrices $\bo{A}$, $\bo{B}$ and $\bo{C}$ is based on the residuals. For the spatial sign score, obvious estimates are
\begin{eqnarray*}
\hat{\bo{A}} = \ave \left [ \bo{A}(\bo{y}_i-\hat{\bs{\beta}}'\bo{x}_{i})\right ], \hspace{5pt}
\hat{\bo{B}} = \frac{1}{n} \hat{\bo{T}}'\hat{\bo{T}} \hbox{ and }
\hat{\bo{C}} = \frac{1}{k} \hat{\bo{T}}'(\bo{Z}\bo{Z}'-\bo{I}_n)\hat{\bo{T}}
\end{eqnarray*}
where $k=\bo{1}_n'\left( \bo{Z}\bo{Z}'-\bo{I}_n \right )\bo{1}_n$. For spatial ranks one could use $\hat{\bo{A}}= \ave [ \bo{A} ( \bo{y}_i -\bo{y}_j -\hat{\bs{\theta}}_{sr}  )  ]$
in which $(\bo{Z}\bo{Z}')_{ij}=0$ and observation $i$ belongs to sample $r$ and $j$ to sample $s$.
An alternative and simpler route is to estimate (\ref{eqn:covmat_of_pairwise_diff}) by a similar estimate as in the one-sample case: $$\hat{\bo{A}}^{-1} \left ( n \hat{\bo{T}}'{\bo{W}^*_{ij}} \bo{Z} \bo{Z}' {\bo{W}^*_{ij}} \hat{\bo{T}} \right ) \hat{\bo{A}}^{-1}
\hbox{, where } \bo{W}_{ij}^*=diag(\bo{w}_{ij}^*).$$ This estimate uses only two samples in the estimation of the middle part
and is more reliable when the variances are heterogeneous across samples.

\section{\label{section:simlations}Efficiency Studies}

In this section we focus on the efficiency of two-sample tests.
For earlier efficiency studies of the multivariate one-sample problem with spatial sign and
rank scores we refer to \citet{Larocque2003}, \citet{LarocqueNevalainenOja2007} and \citet{Haataja2009}.

We generated clustered multivariate data from a linear mixed model setting up a trivariate $t_{\nu}$-distribution.
Full details of the model, designs and cluster size configurations are given in the supplemental file.

The performance of six two-sample tests, Hotelling's $T^2$, spatial sign test, spatial rank test, and their weighted versions was investigated. The weights optimal for the classical Hotelling's $T^2$ test were used for all the weighted tests. Practical experience has shown that the three optimal weight matrices of the tests are very similar, and the rationale for choosing these weights lies in their appealing ease of computation. The tests were studied under three frequently encountered designs (section \ref{subsection:severalsamplestesting}) for different values of the intracluster correlation $\rho$.

\subsection{Asymptotic relative efficiency}

\begin{figure}[h!]
\includegraphics[width=\textwidth]{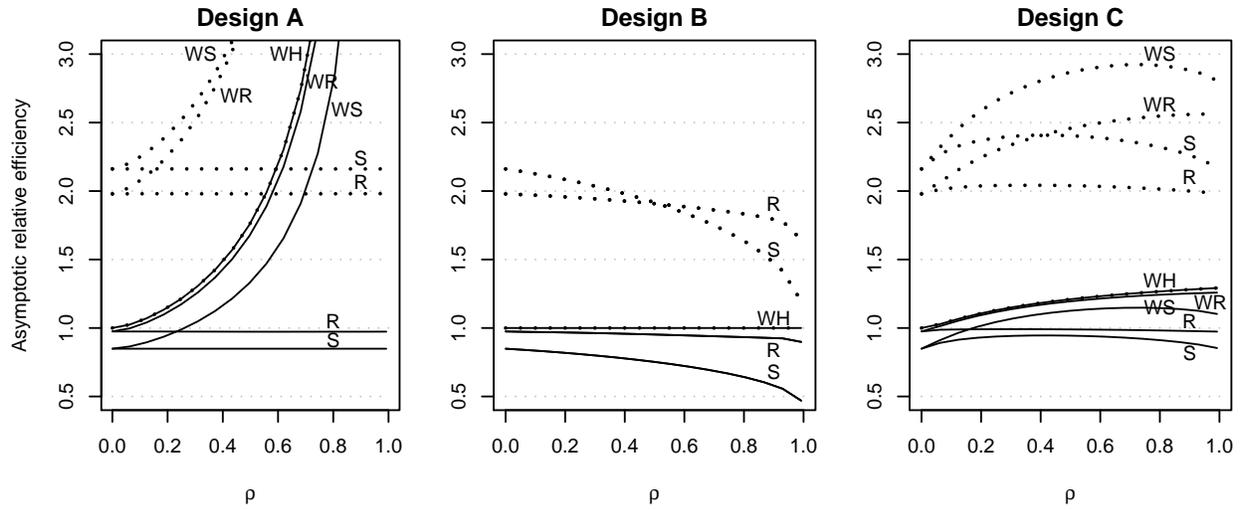}
\caption{\label{figure:are}AREs of different tests relative to Hotelling's $T^2$ test for different designs under
trivariate $t_3$ (dotted lines) and normal distributions (solid lines). The tests are spatial sign tests (S), spatial rank tests (R),
weighted Hotelling's $T^2$ tests (WH), weighted spatial sign tests (WS) and weighted spatial rank tests (WR).}
\end{figure}

Asymptotic relative efficiencies (ARE), using the unweighted Hotelling's $T^2$ as the benchmark test, for the three different designs are shown in Figure \ref{figure:are}. At $\rho=0$, the tests inherit the efficiencies from the i.i.d. case. Hotelling's $T^2$ test is the optimal for the multinormal distribution, but the spatial sign test is the best for the $t_3$-distribution. The ARE of the weighted Hotelling's $T^2$ relative
to the unweighted Hotelling's $T^2$ does not depend on the degrees of freedom. The spatial rank test has a good ARE for both distributions.  The behavior of the tests is remarkably different from design to another when $\rho>0$.

In Design A, the AREs of the unweighted tests do not depend on $\rho$, because $\bo{D}_C$ is here a zero matrix. Weighted tests behave gorgeously, however. Optimal weighting assigns large weights to observations in clusters with both groups present, and less weight to clusters with members only from one group. As $\rho \rightarrow 1$, these within-cluster comparisons tend to receive all the weight, because the treatment effect can be recovered most accurately from them. The unweighted tests still suffer from the error in between-clusters comparisons, and thus the AREs of the weighted tests are dramatically better. Design B is an example of a design where the efficiency cannot be improved by weighting. Thus, the efficiencies of the unweighted and weighted tests overlap. For both spatial sign and rank tests the AREs decrease as a function of $\rho$. The AREs in Design C are similar to the ones in the one-sample case \citep[\emph{e.g}]{LarocqueNevalainenOja2007}. AREs of the unweighted tests at $\rho=0$ and $\rho=1$ are identical since a cluster becomes a singleton observation at $\rho=1$. In between, spatial sign and rank tests suffer less from intracluster correlation than Hotelling's $T^2$. Notable improvements can be obtained by weighting.

\subsection{Simulations}

Results on the empirical size and power are presented in Table \ref{table:simulation results1} for $d=30$. All six tests generally maintain their nominal size well with a few exceptions. Hotelling's $T^2$ and its weighted version are conservative for the $t_3$ error distribution.
The distribution has very heavy tails, and in such a setting the convergence of the moment-based test to its limiting distribution
is slower. All unweighted tests seem conservative in Design C. This suggests that designs without within-cluster comparisons between treatments need a larger sample size for a good $\chi^2_3$-approximation. Interestingly, corresponding weighted tests are liberal
for the same design in particular with larger values of $\rho$. Of course, with $\rho=0.05$ the unweighted and the weighted test are almost the same.
It is worth noting that at $d=14$ (supplemental Table 3), spatial sign and rank tests still maintain their size fairly well in Design A, regardless of the values of $\nu$ and $\rho$, but not as generally in Designs B and C. Curiously, the direction from which the tests converge to meet their target level seems to depend on all parameters of the configuration: score, design, weights, error distribution and intracluster correlation. At $d=60$ the levels overall have improved substantially (supplemental Table 4).

Spatial sign and rank tests have a good power among all the studied error distributions and designs, whereas Hotelling's $T^2$ is the best unweighted
test in the normal case, but has almost no power for the $t_3$-distribution. This problem cannot be solved by selecting a different design, or by weighting. As for the AREs, use of weights enhance the power of the tests for large values of $\rho$: modestly in Design A and more effectively in Design C. No gains of power can be obtained by weighting in Design B.

\section{\label{section:conclusions}Concluding Remarks}

It is commonly acknowledged that in the one-sample case ignoring a positive intracluster correlation leads to
too liberal analyses because of underestimation of the variance. However, our variance breakdown and the supplemental example demonstrate that the opposite can also take place in a multi-treatments study. Therefore, the analysis may either be too liberal or too conservative depending on the design of the study.

The proposed score-based procedures are not necessarily affine invariant and equivariant. Affine invariant
and equivariant versions of spatial sign and rank methods can be developed by using the well-established
transformation-retransformation techniques with an inner standardization \citep{ChakrabortyChaudhuri1996,ChakrabortyChaudhuri1998,ChakrabortyChaudhuri1998b}.  \citet{Larocque2003} and \citet{NevalainenLarocqueOja2007}
present ideas how to do so with clustered data for the testing and estimation problems, respectively. These computationally intensive modifications of the transformation-retransformation procedures can be used as general tools to achieve affine invariance/equivariance with these type of data.

This paper provided a general treatment of score-based testing and estimation methods for multivariate location problems
with clustered data. Asymptotic results were confirmed with simulation studies, which also clearly demonstrate the
gains obtained by the use of scores and weights. In future research we intend to work with multilevel or hierarchical data, and on regression problems.

\section*{Acknowledgements}
Constructive comments by the Associate Editor and two anonymous referees greatly improved the paper.
This research was supported by the Academy of Finland, The Finnish Foundation of Cardiovascular Research and the Competitive Research Funding of the Pirkanmaa Hospital District. The research work of Denis Larocque was supported by the Natural Sciences and Engineering Research Council of Canada.

\section{\label{section:supplemental}Supplemental Materials}
\begin{description}
\item[Supplemental tables:] Details for the efficiency studies and simulation results with $d=14$ and $d=60$. (.pdf)
\item[R-functions used for the simulations:] Functions that generate data from the model of section \ref{section:simlations},
perform the tests, and collect the results into external files. (.zip)
\item[Examples:] Description and analysis of two example data sets. (.pdf)
\item[R-functions used for the examples:] Functions that perform the analyses of the example data sets. (.zip)
\end{description}

\begin{table}[h!]
\begin{center}
\begin{scriptsize}
\caption{\label{table:simulation results1}Empirical size and power of the six two-sample tests under trivariate $t$-distributions with 30 clusters.}
\begin{tabular}{lccccccccccccrr}
  \hline
         & &  \multicolumn{13}{c}{$\rho=0.05$}\\
  \hline
         & &  \multicolumn{6}{c}{$\bs{\Delta}=\bo{0}$} & & \multicolumn{6}{c}{$\bs{\Delta}=\bs{\Delta}_0/\sqrt{N}$}\\
  \cline{3-8}\cline{10-15}
  Design & $\nu$& H & S & R & WH & WS & WR & & H & S & R & WH & WS & WR \\
  \hline
  A    & 3             & 0.044 & 0.050 & 0.051 & 0.044 & 0.051 & 0.050 & & 0.084 & 0.779 & 0.553 & 0.085 & 0.779 & 0.555\\
       & 10            & 0.047 & 0.047 & 0.047 & 0.046 & 0.047 & 0.048 & & 0.560 & 0.798 & 0.755 & 0.561 & 0.799 & 0.757\\
       & $\infty$      & 0.048 & 0.046 & 0.045 & 0.047 & 0.045 & 0.046 & & 0.873 & 0.808 & 0.864 & 0.875 & 0.811 & 0.867\\
  \hline
  B        & 3         & 0.043 & 0.050 & 0.050 & 0.043 & 0.050 & 0.050 & & 0.082 & 0.792 & 0.564 & 0.082 & 0.792 & 0.564\\
           & 10        & 0.048 & 0.048 & 0.050 & 0.048 & 0.048 & 0.050 & & 0.571 & 0.810 & 0.767 & 0.571 & 0.810 & 0.767\\
           & $\infty$  & 0.048 & 0.050 & 0.049 & 0.048 & 0.050 & 0.049 & & 0.886 & 0.825 & 0.880 & 0.886 & 0.825 & 0.880 \\
  \hline
  C          & 3  &      0.039 & 0.045 & 0.043 & 0.042 & 0.049 & 0.047 & & 0.060 & 0.551 & 0.356 & 0.061 & 0.567 & 0.374\\
             & 10 &      0.041 & 0.044 & 0.042 & 0.046 & 0.047 & 0.046 & & 0.371 & 0.575 & 0.520 & 0.385 & 0.593 & 0.541\\
             & $\infty$ & 0.044 & 0.046 & 0.045 & 0.047 & 0.049 & 0.048 & & 0.634 & 0.580 & 0.625 & 0.652 & 0.600 & 0.644\\
  \hline
         & &  \multicolumn{13}{c}{$\rho=0.2$}\\
  \hline
  A    & 3             & 0.044 & 0.050 & 0.052 & 0.043 & 0.051 & 0.049 & & 0.084 & 0.793 & 0.575 & 0.083 & 0.817 & 0.602\\
       & 10            & 0.049 & 0.049 & 0.048 & 0.048 & 0.048 & 0.048 & & 0.579 & 0.819 & 0.778 & 0.598 & 0.840 & 0.805\\
       & $\infty$      & 0.047 & 0.049 & 0.047 & 0.048 & 0.046 & 0.047 & & 0.891 & 0.826 & 0.883 & 0.917 & 0.852 & 0.908\\
  \hline
  B        & 3         & 0.044 & 0.050 & 0.049 & 0.044 & 0.050 & 0.049 & & 0.086 & 0.844 & 0.630 & 0.086 & 0.844 & 0.630\\
           & 10        & 0.049 & 0.050 & 0.050 & 0.049 & 0.050 & 0.050 & & 0.631 & 0.864 & 0.834 & 0.631 & 0.864 & 0.834\\
           & $\infty$  & 0.050 & 0.049 & 0.050 & 0.050 & 0.049 & 0.050 & & 0.936 & 0.877 & 0.928 & 0.936 & 0.877 & 0.928\\
  \hline
  C          & 3  &       0.038 & 0.044 & 0.041 & 0.053 & 0.056 & 0.057 & & 0.052 & 0.345 & 0.216 & 0.061 & 0.405 & 0.263\\
             & 10 &       0.041 & 0.044 & 0.042 & 0.057 & 0.057 & 0.057 & & 0.233 & 0.369 & 0.320 & 0.282 & 0.423 & 0.378\\
             & $\infty$ & 0.044 & 0.044 & 0.044 & 0.059 & 0.058 & 0.057 & & 0.384 & 0.368 & 0.389 & 0.451 & 0.422 & 0.449\\
  \hline
         & &  \multicolumn{13}{c}{$\rho=0.4$}\\
  \hline
  A    & 3             & 0.044 & 0.050 & 0.052 & 0.046 & 0.052 & 0.053 & & 0.086 & 0.817 & 0.607 & 0.087 & 0.888 & 0.693\\
       & 10            & 0.049 & 0.050 & 0.049 & 0.053 & 0.051 & 0.053 & & 0.607 & 0.844 & 0.808 & 0.680 & 0.915 & 0.894\\
       & $\infty$      & 0.048 & 0.051 & 0.048 & 0.054 & 0.050 & 0.053 & & 0.914 & 0.850 & 0.903 & 0.970 & 0.923 & 0.965\\
  \hline
  B        & 3         & 0.045 & 0.053 & 0.052 & 0.045 & 0.053 & 0.052 & & 0.092 & 0.914 & 0.739 & 0.092 & 0.914 & 0.739\\
           & 10        & 0.057 & 0.054 & 0.057 & 0.057 & 0.054 & 0.057 & & 0.724 & 0.935 & 0.918 & 0.724 & 0.935 & 0.918\\
           & $\infty$  & 0.058 & 0.052 & 0.057 & 0.058 & 0.052 & 0.057 & & 0.980 & 0.945 & 0.977 & 0.980 & 0.945 & 0.977\\
  \hline
  C          & 3  &  0.037 & 0.043 & 0.042 & 0.062 & 0.064 & 0.066 & & 0.046 & 0.231 & 0.148 & 0.067 & 0.304 & 0.204     \\
             & 10 &  0.041 & 0.044 & 0.043 & 0.068 & 0.065 & 0.066 & & 0.157 & 0.249 & 0.213 & 0.223 & 0.319 & 0.285 \\
             & $\infty$ & 0.041 & 0.043 & 0.042 & 0.069 & 0.064 & 0.066 & & 0.249 & 0.248 & 0.253 & 0.336 & 0.315 & 0.333\\
  \hline
  \multicolumn{15}{l}{H = Hotelling's $T^2$ test; S = Spatial sign test; R = Spatial rank test}\\
  \multicolumn{15}{l}{WH = Weighted Hotelling's $T^2$ test; WS = Weighted spatial sign test; WR = Weighted spatial rank test}
\end{tabular}
\end{scriptsize}
\end{center}
\end{table}
\newpage
\newpage

\end{document}